\documentclass[11pt]{article}

\usepackage{authblk}
\usepackage{graphicx}
\usepackage[utf8]{inputenc}
\usepackage[numbers]{natbib}
\usepackage{amsmath}
\usepackage{amsfonts}
\usepackage{amssymb}
\usepackage{caption}
\usepackage{subcaption}
\usepackage[spaces,hyphens]{url}

\title{Stability Analysis of Fractional Difference Equations with Delay}
\author{Divya D. Joshi\\ Department of Physics, Rashtrasant Tukadoji Maharaj Nagpur University, Amravati Rd, Nagpur-
440033, Maharashtra, India
	\and Sachin Bhalekar\\ School of Mathematics and Statistics, University of Hyderabad, Hyderabad-500046, Telangana, India
	\and Prashant M. Gade\\ Department of Physics, Rashtrasant Tukadoji Maharaj Nagpur University, Amravati Rd, Nagpur-
	440033, Maharashtra, India}

\date{\today}

\newtheorem{Lem}{Lemma}[section]

\newtheorem{Def}{Definition}[section]
\newtheorem{Pro}{Property}[section]

\begin{document}
\maketitle



\begin{abstract}
	Long-term memory is a feature observed in systems ranging from neural networks to epidemiological models. The memory in such systems is usually modeled by the time delay. Furthermore, the nonlocal operators, such as the "fractional order difference" can also have a long-time memory. Therefore, the fractional difference equations with delay are an appropriate model in a range of systems. Even so, there are not many detailed studies available related to the stability analysis of fractional order systems with delay. In this work, we derive the stability conditions for linear fractional difference equations with a delay term $\tau$. We have given detailed stability analysis for the cases $\tau=1$ and $\tau=2$. The results are extended to nonlinear maps.
\end{abstract}



	\section{Introduction}
The Discovery of calculus by Newton and Leibnitz revolutionized natural sciences. Since then, differential equations have been a substantial tool for modeling various phenomena for mathematicians, physicists, and engineers. This tool is inadequate for systems that are non-local in nature. They have been modeled by fractional differential/difference equations. They show hereditary properties and have memory. For example, Bagley and Torvik have studied the theoretical model of viscoelastic materials defined by fractional derivatives for polymer solutions and polymers solids without crosslinking \cite{bagley1983theoretical,bagley1986fractional}. 
Fractional differential equations have been also used to investigate epidemiological topics such as the transmission of COVID-19 \cite{baleanu2020fractional, zhang2020dynamics, chatterjee2021fractional}, Ebola virus \cite{srivastava2020efficient} etc.
Besides, fractional calculus is found to be a great foundation for modeling non-conventional fractal and non-local media and fractional operators are useful at explaining complex long-memory and multiscale phenomena in materials \cite{failla2020advanced}. According to Torvik and Bagley, most real-life phenomena
can be modeled by fractional calculus \cite{torvik1984appearance}.

Another set of equations that can be used to model memory-based systems is delay differential equations. They have found applications in a wide range of different fields. Mackey and Glass presented delay differential equations to model complex dynamics in physiological control systems \cite{glass1988clocks}. A differential delay equation has been presented to model oscillations in laboratory blowfly population \cite{perez}. For an optical bistable resonator, Ikeda and Matsumoto showed that dynamics can be reproduced with a single bistable resonator \cite{ikeda1987high}. The differential equation with delay is effectively an infinite dimensional system and Lu and He showed that difference-differential equations with a single variable can show chaos for large delay \cite{lu1996chaotic}. Delays and their impact on the stability of ecological systems is a widely debated topic \cite{eurich2005distributed}. Delay differential equations are widely used in topics as diverse as chemical	kinetics \cite{roussel1996use} to population dynamics of agricultural pests \cite{schley2003delay}.  It has found applications in immunology \cite{rihan2020dynamics}, epidemiology \cite{alzahrani2022repercussions} and neural networks \cite{xu2021bifurcation}. 

Fractional differential equations with delay have also received some attention. The fractional extension of the Ikeda system is studied in \cite{jun2006chaotic}. Fractional extension of logistic and Chen equations have also been investigated \cite{wang2008chaos,daftardar2012dynamics}. On a theoretical side, the stability conditions for fractional differential equations with delay are also studied \cite{li2011survey,bhalekar2016stability,bhalekar2019analysis,bhalekar2022can}. Numerical methods to solve these equations are provided in \cite{bhalekar2011predictor,daftardar2015solving}. 

In the qualitative theory of dynamical systems, studies in differential equations and difference equations have gone hand-in-hand. Period-doubling cascade which is observed in several models and experiments of nonlinear dynamical systems was explained by Feigenbaum using maps \cite{feigenbaum1978quantitative}. Almost all routes to chaos have been studied using maps \cite{ott2002chaos}. Several techniques to control chaos are applicable in both differential equations as well as maps \cite{scholl2008handbook}. Thus the analysis of discrete-time systems can render useful insights in continuous time systems. Hence we study fractional difference equations with delay in this work with a focus on the basic problem of stability of the linear system.

Difference equations with delay have been studied since long. The well-investigated 2D maps such as H{\'e}non and Lozi map can be viewed as one-dimensional maps with delay \cite{henon,lozi}. A logistic map with delayed feedback is studied in \cite{buchner2000logistic}. Systems with delay can be viewed as spatiotemporal systems and we can study those from viewpoint of phase transitions. Such studies are carried out in \cite{ashwini,dutta}. 
In integer order difference equations, the stability is determined by whether all eigenvalues of Jacobian are inside a unit	circle in a complex plane. Integer order difference equations with delay can be treated as higher dimensional difference equations and if the eigenvalues are inside the unit circle in a complex plane, it is stable.  We define fractional delay difference equations in a manner analogous to fractional delay differential equations. We find that the criterion is more complicated in this case.
The stability curve need not be a simple curve. The stability region can be a union of distinct disjoint regions. Nonetheless, the analytic criterion can be given in a closed functional form.

Wu and Baleanu have numerically studied fractional-order logistic maps with delay in \cite{wu2015discrete}. Similar studies are carried out in the context of economic model \cite{tarasova2017logistic}. Local stability of the fractional ordered delay Mackey-Glass equation has been analyzed in \cite{el2016fractional} with a different discretization process. The importance of studying delayed systems from the viewpoint of applications as well as control cannot be overemphasized. We extend these studies to exact results for linear fractional order systems and discuss an extension to nonlinear maps.

Certain	inequalities have been derived for fractional order delay difference equations \cite{alzabut2018generalized,chen2021finite,wang2021finite,du2019finite,du2022new}. However, detailed bifurcation analysis of possible routes to instability has not been carried out to the best of our knowledge. In this work, we present a detailed bifurcation analysis for such a system. Bifurcation analysis	gets more complex for a larger delay and detailed bifurcation analysis can be carried out only for a small delay. However, our results are generic and allow us to study a linear system with any delay. The analysis is also generic and can be easily generalized to the system with multiple delays.

\section{Preliminaries}

\begin{Def}(see \cite{mozyrska2015transform})
	The Z-transform of a sequence $ \{y(n)\}_{n=0}^\infty $ is a complex function given by
	\begin{equation*}
		Y(z)=Z[y](z)=\sum_{k=0}^{\infty} y(k) z^{-k},
	\end{equation*}
	where $z \in \mathbb{C}$ is a complex number 
	for which the series converges absolutely. 
\end{Def}
\begin{Def}(see \cite{ferreira2011fractional,bastos2011discrete})
	Let $ h > 0 ,\; a \in \mathbb{R}$ and 
	$ (h\mathbb{N})_a = \{ a, a+h, a+2h, \ldots\} $.
	For a function $x : (h\mathbb{N})_a \rightarrow  \mathbb{C}$, 
	the forward h-difference operator is defined as 
	$$
	(\Delta_h x)(t)=\frac{x(t+h)- x(t)}{h},$$
	where t	$ \in (h\mathbb{N})_a $.
\end{Def}
Throughout this article, we take a = 0 and 
h = 1. We write $\Delta$ for $\Delta_1 $.
Now, we generalize the fractional order operators 
defined in \cite{mozyrska2015transform,ferreira2011fractional,bastos2011discrete}.
\begin{Def}
	For a function  $x : (h\mathbb{N})_a \rightarrow  \mathbb{C}$ 
	the fractional h-sum of order 
	$\alpha = u +\iota v \in \mathbb{C}, u>0 $ is given by
	\begin{equation*}
		(_{a}\Delta_h^{-\alpha}x)(t) 
		= \frac{h^\alpha}{\Gamma(\alpha)}\sum_{s=0}^{n}\frac{\Gamma(\alpha+n-s)}{\Gamma(n-s+1)} x(a+sh),\\
	\end{equation*}
	where,	$t=a+(\alpha+n)h, \; n \in \mathbb{N_\circ}$.
\end{Def}
For $h=1$ and $a=0$	, we have
\begin{eqnarray*}
	(\Delta^{-\alpha}x)(t) 	&=&\frac{1}{\Gamma(\alpha)}\sum_{s=0}^{n}\frac{\Gamma(\alpha+n-s)}{\Gamma(n-s+1)}x(s)\\
	&=&\sum_{s=0}^{n}
	\left(
	\begin{array}{c}
		n-s+\alpha-1\\
		n-s\\
	\end{array}
	\right)
	x(s).
\end{eqnarray*}                         
Here, we used the generalized binomial coefficient
\begin{equation*}
	\left(
	\begin{array}{c}
		\mu \\
		\eta\\
	\end{array}
	\right)
	=\frac{\Gamma(\mu+1)}{\Gamma(\eta+1)\Gamma(\mu-\eta+1)},\\
	\; \mu ,\eta \in \mathbb{C}, \; \text{Re}(\mu)>0,\;\text{and Re}(\eta)>0.
\end{equation*}
If $n$ $\in$ $\mathbb{N_\circ}$ then
\begin{eqnarray*}
	\left(
	\begin{array}{c}
		\mu \\
		n
	\end{array}
	\right)
	=\frac{(\mu + 1)}{n!\Gamma(\mu-n+1)}
	=\frac{\mu(\mu-1)\ldots(\mu-n-1)}{n!}.
\end{eqnarray*}
\begin{Def}
	For $n \in \mathbb{N_\circ}$ and $\alpha=u+\iota v \in \mathbb{C}, u>0,$ we define
	\begin{eqnarray*}
		\tilde{\phi}_{\alpha}(n)=
		\left(
		\begin{array}{c}
			n+\alpha-1\\
			n\\
		\end{array}
		\right)
		=(-1)^n
		\left(
		\begin{array}{c}
			-\alpha\\
			n
		\end{array}
		\right).
	\end{eqnarray*}
\end{Def}
\textbf{Note}: The convolution $\tilde{\phi}_{\alpha}*x$ of the sequences $\tilde{\phi}_{\alpha}$ and $x$ is defined as
\begin{equation*}
	\left(\tilde{\phi}_{\alpha}*x\right)(n)=\sum_{s=0}^{n}\tilde{\phi}_{\alpha}(n-s)x(s)
\end{equation*}
\begin{equation*}
	\therefore (\Delta^{-\alpha}x)(n)=(\tilde{\phi}_{\alpha}*x)(n).\\
\end{equation*}
\begin{eqnarray*}
	\therefore Z(\Delta^{-\alpha}x)(n)=Z\left(\tilde{\phi}(n)\right)Z(x(n))\\
	=(1-z^{-1})^{-\alpha}X(z),	
\end{eqnarray*}
where $X$ is $Z$ transform of $x$.
\begin{Pro}(see \cite{petalez})
	The time-shifting property shows how a change in the discrete function's time domain alters the Z-domain.
\end{Pro}
\begin{equation*}
	Z[x(k-n)]=z^{-n}X(z)
\end{equation*}
\textbf{Proof}:
From Definition 2.1 we have,
\begin{equation*}
	X(z)=\sum_{k=0}^{\infty}x(k)z^{-k}.
\end{equation*}
Consider $k-n=m$ i.e., $k=m+n$. Thus, we write the z-transform equation as 
\begin{equation*}
	Z[x(k-n)]=\sum_{k=0}^{\infty}x(k-n)z^{-k}\\
	=\sum_{m=0}^{\infty}x(m)z^{-(m+n)}\\
	=\sum_{m=0}^{\infty}x(m)z^{-m}z^{-n}\\
	=z^{-n}\sum_{m=0}^{\infty}x(m)z^{-m}\\
	=z^{-n}X(z).
\end{equation*}
\begin{Lem}
	For $\alpha \in \mathbb{C},\; \text{Re}(\alpha)>0$,
	\begin{equation*}
		Z(\tilde{\phi}_{\alpha}(t))=\frac{1}{(1-z^{-1})^{\alpha}}.
	\end{equation*}
\end{Lem}
\textbf{Proof}: We have, 
\begin{eqnarray*}
	Z(\tilde{\phi}_{\alpha}(t))&=&\sum_{j=0}^{\infty}\tilde{\phi}_{\alpha}(j)z^{-j}\\
	&=&\sum_{j=0}^{\infty}\left(
	\begin{array}{c}
		j+\alpha-1\\
		j
	\end{array}
	\right)z^{-j}\\
	&=&\sum_{j=0}^{\infty}(-1)^{j}\left(
	\begin{array}{c}
		-\alpha\\
		j
	\end{array}
	\right)z^{-j}\\
	&=&(1-z^{-1})^{-\alpha}.
\end{eqnarray*}
by using Newton's generalization of the Binomial Theorem 
\cite{niven1969formal,link1}.
\section{Model}
Consider the fractional order linear difference equation
\begin{equation}
	x(t)=x_0+\sum_{j=0}^{t-1}\frac{\Gamma(t-j+\alpha-1)}{\Gamma(\alpha)\Gamma(t-j)}(a-1)x(j), \label{eqn1}
\end{equation}
where, $x_{0}=x(0)$.
In this paper, we study the stability analysis of a fractional difference equation (\ref{eqn1}) with a delay term. 
\subsection*{Modeling Equation}
Introducing the delay term in the equation as $bx(t-\tau)$, we get
\begin{equation}
	x(t)=x_{0}+\sum_{j=0}^{t-1}\left(\frac{\Gamma(t-j+\alpha-1)}{\Gamma(\alpha)\Gamma(t-j)}\left((a-1)x(j)+b x(j-\tau)\right)\right), \label{eqn2}
\end{equation}
where, $b \in \mathbb{R}$ and $a \in \mathbb{C}$.
\begin{equation*}
	\therefore	x(t+1)=x_0+(a-1)(\tilde{\phi}_{\alpha}*x)(t)+b (\tilde{\phi}_{\alpha}*x)(t-\tau).
\end{equation*}
Taking Z-transform on both sides, we get
\begin{eqnarray}
	zX(z)-zx_0 \nonumber	&=&\frac{x_0}{1-z^{-1}}+\frac{(a-1)}{(1-z^{-1})^{\alpha}}X(z)+\frac{b z^{-\tau}}{(1-z^{-1})^{\alpha}}X(z). \nonumber\\	
	\therefore	X(z)\left(z-\frac{(a-1)}{(1-z^{-1})^{\alpha}}+\frac{b z^{-\tau}}{(1-z^{-1})^{\alpha}}\right)&=&\frac{x_0}{1-z^{-1}}+zx_0,\label{eqn3}
\end{eqnarray} where $|z|<1$.
\subsection*{Characteristic Equation}
From (\ref{eqn3}), the characteristic equation of (\ref{eqn2}) is
\begin{equation}
	\left(z(1-z^{-1})^{\alpha}-(a-1)-b z^{-\tau}\right)=0, \label{eqn4}
\end{equation}
where the condition $|z|<1$ should be satisfied. 
Putting $z=e^{(it)}$ in equation (\ref{eqn4}), we get,
\begin{eqnarray}
	e^{(it)}(1-e^{-it})^{\alpha}-be^{-it\tau}-(a-1)&=&0,  \nonumber \\
	\text{i.e. }
	e^{(it)}(1-e^{-it})^{\alpha}-be^{-it\tau}+1&=& a.\label{eqn 5}
\end{eqnarray}
\subsection*{Matrix Representation}
Equation (\ref{eqn2}) can be represented equivalently as the following system
\begin{eqnarray}
	x(t)&=&x_0+\sum_{j=0}^{t-1}\frac{\Gamma(t-j+\alpha-1)}{\Gamma(\alpha)\Gamma(t-j)}\left((a-1)x(j)+b y(j)\right),\; y(t)=x(t-1) \text{ for } \tau=1.\label{eqn 6}\\
	\therefore\; x(t+1)&=&x_0+(a-1)(\tilde{\phi}_{\alpha}*x)(t)+b (\tilde{\phi}_{\alpha}*y)(t). \nonumber
\end{eqnarray}
Taking z-transform, we get
\begin{equation*}
	zX(z)-zx_0	=\frac{x_0}{1-z^{-1}}+\frac{(a-1)}{(1-z^{-1})^{\alpha}}X(z)+\frac{b}{(1-z^{-1})^{\alpha}}, \; Y(z)=\frac{X(z)}{z}+x(-1).
\end{equation*}
\begin{equation*}
	\therefore	(z(1-z^{-1})^{\alpha}-(a-1))X(z)-bY(z) = x_0(1-z^{-1})^{\alpha}(z-(1-z^{-1})^{-1}),
\end{equation*}
\begin{eqnarray*}
	zY(z)-X(z)=zx(-1).\\
	\therefore	\begin{bmatrix}
		z(1-z^{-1})^{\alpha}-(a-1) & -b\\
		-1 & z
	\end{bmatrix}
	\begin{bmatrix}
		X(z)\\
		Y(z)
	\end{bmatrix}
	=0.\\
	\therefore	\begin{vmatrix}
		z(1-z^{-1})^{\alpha}-(a-1) & -b\\
		-1 & z
	\end{vmatrix}
	=0.\\
	\text{i.e. }	z(1-z^{-1})^{\alpha}-b z^{-1}+1=a.
\end{eqnarray*}
This is the equation for the model with $\tau=1$.
Similarly, for $\tau=2$ we get,
\begin{eqnarray*}
	\begin{vmatrix}
		z(1-z^{-1})^{\alpha}-(a-1) & 0 & -b\\
		-1 & (z-1)+1 & 0\\
		0 & -1 & (z-1)+1
	\end{vmatrix}
	=0.\\
\end{eqnarray*}
We can generalize and write $\tau+1$ 
dimensional determinant for delay $\tau$  as follows,

\begin{eqnarray*}
	\begin{vmatrix}
		z(1-z^{-1})^{\alpha}-(a-1) & 0 & 0 & $\ldots$ & -b\\
		-1 & (z-1)+1 & 0 & \ldots & 0 \\
		0 & -1 & (z-1)+1  & \ldots & 0 \\
		\vdots & \ddots & \ddots & \ldots & \vdots \\
		0 & \ldots & -1  & (z-1)+1 & 0\\
		0 & 0 & \ldots  & -1 & (z-1)+1
	\end{vmatrix}
	=0.\\
\end{eqnarray*}

\subsection*{Boundary curve}
We give the parametric representation of the boundary condition (\ref{eqn 5}) by,
\begin{equation}
	\begin{aligned}
		\gamma(t)= & \{2^{\alpha}\left(\sin\left(\frac{t}{2}\right)\right)^{\alpha}\cos\left(\frac{\alpha \pi}{2}+\left(1-\frac{\alpha}{2}\right)\right)-b\cos(\tau t)+1, \\
		& 2^{\alpha}\left(\sin\left(\frac{t}{2}\right)\right)^{\alpha}\sin\left(\frac{\alpha \pi}{2}+\left(1-\frac{\alpha}{2}\right)\right)+b\sin(\tau t)\} ,
	\end{aligned}\label{eqn 7}
\end{equation}
for $t \in [0,2\pi]$ in complex plane.
If the complex number $a$ lies inside this anticlockwise oriented simple closed curve $\gamma(t)$ then
the system (\ref{eqn2}) will be asymptotically stable.

\subsection{Stability Analysis with $\tau=1$} \label{sec3}
\begin{figure*}[h]
	\centering
	\includegraphics[scale=0.7]{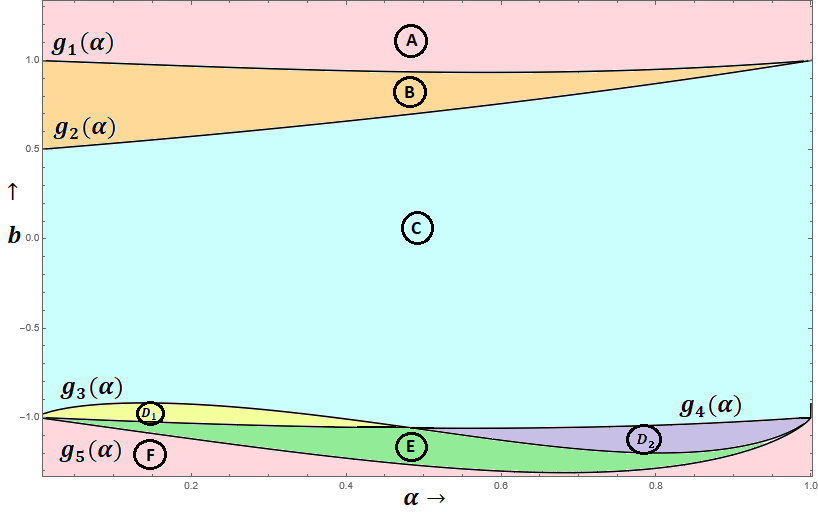}
	\caption{bifurcation regions in $\alpha b$-plane for
		$\tau=1$.}
	\label{fig1}
\end{figure*}

\begin{figure}
	\centering\subfloat[Region $A$]{%
		\includegraphics[scale=0.24]{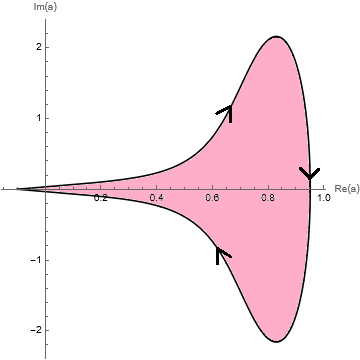}
	}
	\hspace{0.5cm}
	\centering\subfloat[Region $B$]{%
		\includegraphics[scale=0.2]{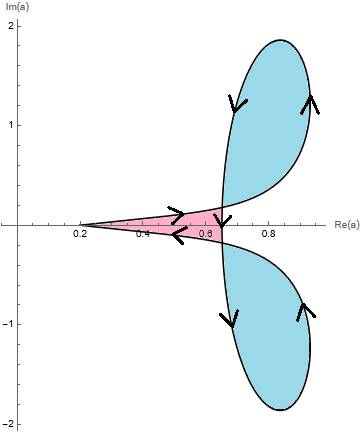}
	}
	\hspace{0.5cm}
	\centering\subfloat[Region $C$]{%
		\includegraphics[scale=0.22]{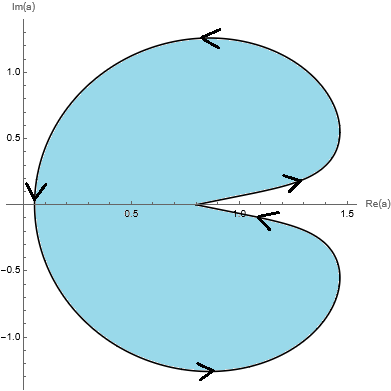}
	}
	\hspace{0.5cm}
	\centering\subfloat[Region $D_{1}$]{%
		\includegraphics[scale=0.21]{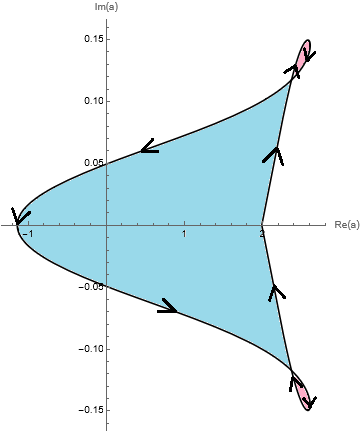}
	}\\
	\centering\subfloat[Region $D_{2}$]{%
		\includegraphics[scale=0.21]{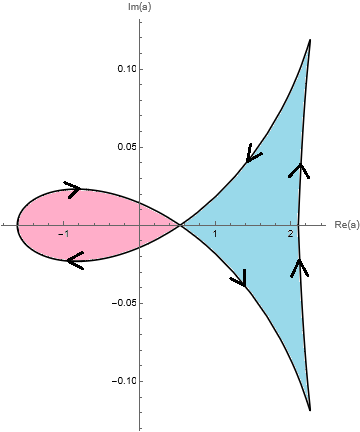}
	}
	\hspace{0.5cm}
	\centering\subfloat[Region $E$]{%
		\includegraphics[scale=0.26]{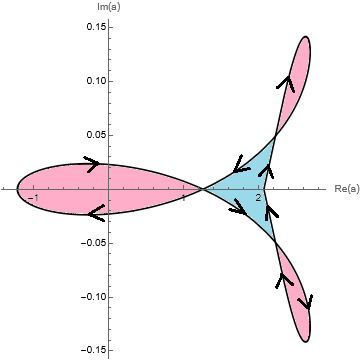}
	}
	\hspace{0.5cm}
	\centering\subfloat[Region $F$]{%
		\includegraphics[scale=0.28]{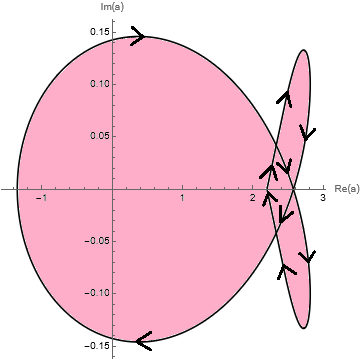}
	}
	\hspace{0.5cm}
	\centering\subfloat[Curve $b=g_{1}(\alpha)$]{%
		\includegraphics[scale=0.26]{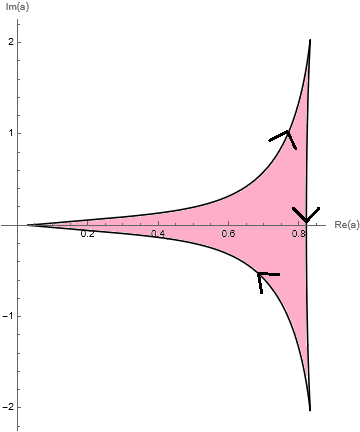}
	}\\
	\centering\subfloat[Curve $b=g_{2}(\alpha)$]{%
		\includegraphics[scale=0.29]{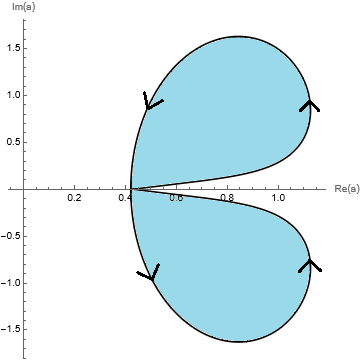}
	}
	\hspace{0.5cm}
	\centering\subfloat[Curve $b=g_{3}(\alpha)$ ($\alpha<\alpha_{*}$)]{%
		\includegraphics[scale=0.26]{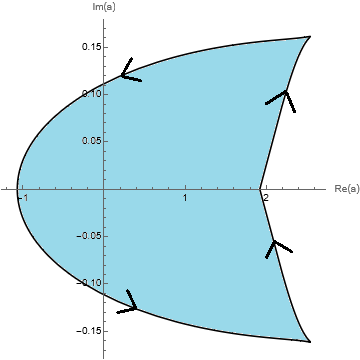}
	}
	\hspace{0.5cm}
	\centering\subfloat[Curve $b=g_{3}(\alpha)$ ($\alpha>\alpha_{*}$)]{%
		\includegraphics[scale=0.23]{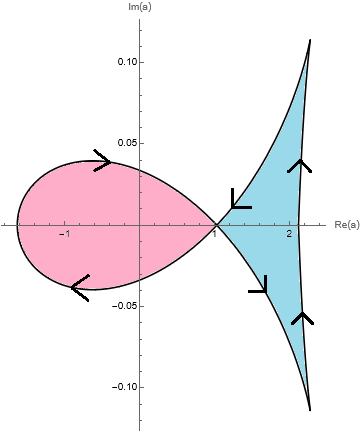}
	}
	\hspace{0.5cm}
	\centering\subfloat[Curve $b=g_{4}(\alpha)$ ($\alpha<\alpha_{*}$)]{%
		\includegraphics[scale=0.3]{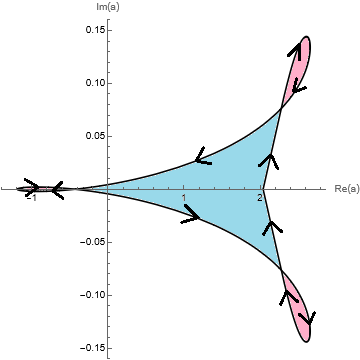}
	}\\
	\centering\subfloat[Curve $b=g_{4}(\alpha)$ ($\alpha>\alpha_{*}$)]{%
		\includegraphics[scale=0.25]{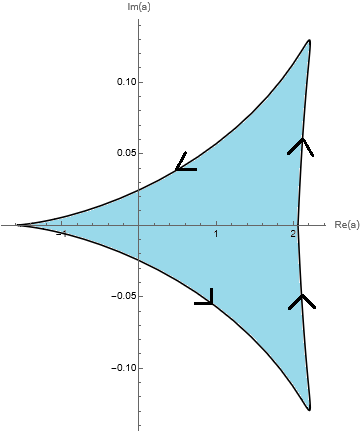}
	}
	\hspace{0.5cm}
	\centering\subfloat[Curve $b=g_{5}(\alpha)$]{%
		\includegraphics[scale=0.29]{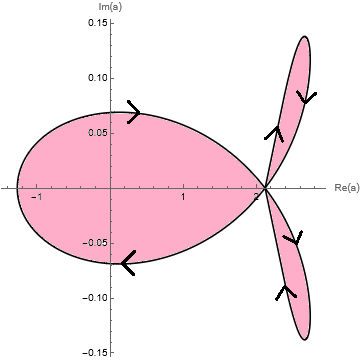}
	}
	\caption{The figures 2(a), 2(b), 2(c), 2(d), 2(e), 2(f), 2(g) are the representative stability diagrams having $b$ and $\alpha$ values lying in the regions $A, B, C, D_1,D_2, E,$ and $F$ shown in the figure (\ref{fig1}) respectively. The figures 2(h), 2(i), 2(j), 2(k), 2(l), 2(m) and 2(n) are the representative stability diagrams having $b$ and $\alpha$ values lying on the branches $b=g_{1}(\alpha)$, $b=g_{2}(\alpha)$, $b=g_{3}(\alpha)$ ($\alpha<\alpha_{*}$),$b=g_{3}(\alpha)$ ($\alpha>\alpha_{*}$), $b=g_{4}(\alpha)$ ($\alpha<\alpha_{*}$), $b=g_{4}(\alpha)$ ($\alpha>\alpha_{*}$), and $b=g_{5}(\alpha)$ respectively. Here the blue color indicates stable regions and the pink color indicates unstable regions.}
	\label{fig2}
\end{figure}

For given $b \in \mathbb{R}$, if the boundary curve $\gamma(t)$ defined by (\ref{eqn 7}) is simple, closed, and anticlockwise oriented (see Figure (\ref{fig2}(c))) then for any region $a \in \mathbb{C}$ in the region bounded by $\gamma(t)$ gives stable solutions to the system (\ref{eqn2}).

However, for some values of $b\in \mathbb{R}$, $\gamma(t)$ is not simple (\ref{fig2}). In this case, we can have multiple regions bounded by $\gamma$, and some of them are negatively oriented. Therefore, the sub-regions with negative orientation give unstable regions. Moreover, there are some values of $b$, for which $\gamma$ is simple but the orientation is clockwise. We get complete instability in these cases.

Our observations are as follows:\\
Self-intersections and cusps are the main routes to the multiple points in $\gamma(t)=\left(\gamma_1(t),\gamma_2(t)\right)$.
At  $b=g_{2}(\alpha)=2^{\alpha-1}$, the multiple point in  $\gamma(t)$ is due to $\gamma(0)=\gamma(\pi)$. We also observe that, there exists $t_{0}\in(0,\pi)$ such that $\gamma(0)=\gamma(t_{0})$ and $\gamma_{2}(t_{0})=0$. Solving $\gamma_{2}(t_{0})=0$ for $b$, we get
\begin{equation}
	b=-2^{\alpha}\csc\left(t_{0}\right)\left(\sin\left(\frac{t_{0}}{2}\right)\right)^{\alpha}\sin\left(t_{0}\left(1-\frac{\alpha}{2}\right)+\frac{\pi\alpha}{2}\right).  \label{t2-b1}
\end{equation}
Substituting this value of  $b$ in (\ref{eqn 7}), we get $\gamma_{1}(t_{0})=\gamma_{1}(0)$.  Solving this equation,  we get $t_{0}=\frac{\pi\left(1-\alpha\right)}{\left(3-\alpha\right)}$.
This gives $b=g_{5}(\alpha)=-2^{\alpha}\left(\cos\left(\frac{\pi}{3-\alpha}\right)\right)^{\alpha}\csc\left(\frac{\pi\left(-1+\alpha\right)}{-3+\alpha}\right)\sin\left(\frac{\pi}{3-\alpha}\right).$
This is a curve in the $\alpha b$-plane, below which the orientation of $\gamma(t)$ becomes clockwise and there is no stable region (cf. Figures \ref{fig1}, \ref{fig2}(g)). For the values of $(\alpha, b)$ above and sufficiently close to the curve $b=g_{5}(\alpha)$, the curve $\gamma(t)$ generates a small stable region and three disconnected unstable regions (see Fig. \ref{fig2}(f)).

At the cuspidal points, $\gamma'(t)$ does not exist. We may expect $x'(t)=y'(t)=0.$ Solving $y'(t)=0$ for $b$, we get
\begin{equation}
	b=-2^{-1+\alpha}\sec\left(t\right)\left(\sin\left(\frac{t}{2}\right)\right)^{-1+\alpha}\left(\left(-1+\alpha\right)\sin\left(\frac{1}{2}\left(t+\pi\alpha-t\alpha\right)\right)+\sin\left(\frac{1}{2}\left(3t+\pi\alpha-t\alpha\right)\right)\right). \label{eqn b1}
\end{equation}
Substituting this $b$ in $x'(t)=0$ we get,
\begin{equation*}
	\cos\left(\frac{1}{2}\left(-t\left(-5+\alpha\right)+\pi\alpha\right)\right)+\left(-1+\alpha\right)\cos\left(\frac{1}{2}\left(-t\left(-3+\alpha\right)+\pi\alpha\right)\right)=0
\end{equation*}
For given $\alpha \in (0,1)$, we can solve this equation for $t$ and substitute in equation (\ref{eqn b1}) to generate three branches viz. $b=g_{1}(\alpha)$, $b=g_{3}(\alpha)$, and $b=g_{4}(\alpha)$ as shown in Figure \ref{fig1}.

All these bifurcation curves $b=g_{k}(\alpha)$, $k=1,2,3,4$ and 5 produce the bifurcation regions $A, B, C, D, E$ and $F$. In region $A$, $b>g_{1}(\alpha)$ and $\gamma(t)$ is simple but clockwise oriented (see Figure \ref{fig2}(a)). The system (\ref{eqn2}) is unstable for any $a \in \mathbb{C}$.

At $b=g_{1}(\alpha)$, there are two cusps is $\gamma(t)$ (see Figure \ref{fig2}(h)). If $g_{2}(\alpha)<b<g_{1}(\alpha)$ ({\it{i.e.,}} region $B$), these two cusps result in multiple regions. There are two stable regions and one unstable region bounded by $\gamma(t)$ (see Figure \ref{fig2}(b)).

At $b=g_{2}(\alpha)$, the unstable region disappears and we have two disjoint stable regions due to $\gamma(0)=\gamma(\pi)=\gamma(2\pi)$ (see Figure \ref{fig2}(i)). In the region $C$, $\gamma(t)$ is an anticlockwise oriented simple curve. Therefore, the stable region is bounded by $\gamma(t)$ (see Figure \ref{fig2}(c)).

Note that the curves $b=g_{3}(\alpha)$ and $b=g_{4}(\alpha)$ intersect at $\alpha_{*}=0.486$. For $0<\alpha<\alpha_{*}$, $g_{4}(\alpha)<g_{3}(\alpha)$ and the region bounded by them is $D_{1}$. For $\alpha_{*}<\alpha<1$, $g_{3}(\alpha)<g_{4}(\alpha)$ and they form the region $D_{2}$ (see Figure \ref{fig1}).

At $b=g_{3}(\alpha)$, there are two cusps in $\gamma(t)$ on right side(see Figure \ref{fig2}(j) and (k)). On the other hand, at $b=g_{4}(\alpha)$, there is only one cusp in $\gamma(t)$ which is on the left side(see Figure \ref{fig2}(l) and (m)).

In the region $D_{1}$, $\gamma(t)$ has one stable region (left side) and two unstable regions (right side)(see Figure \ref{fig2}(d)). There is a stable and an unstable region in $D_{2}$ (see Figure \ref{fig2}(e)).

All the three unstable regions in $D_{1}$ and $D_{2}$ are also there in the region $E$ along with one stable region in the central part of $\gamma(t)$ (see Figure \ref{fig2}(f)).

At $b=g_{5}(\alpha)$, all the self-intersection points discussed above merge, and the stable region disappears completely (see Figure \ref{fig2}(n)). If $b<g_{5}(\alpha)$, all the regions generated by $\gamma(t)$ have clockwise orientation, and hence the system is unstable(see Figure \ref{fig2}(g)).

The quick summary of qualitative changes across the bifurcation curves is as follows: \\
1. For all the parameter values above the curve $b=g_1(\alpha)$ and below the curve $b=g_{5}(\alpha)$  in the $\alpha b$-plane, there is no any stable region bounded by the curve $\gamma(t)$ (see Figures \ref{fig2}(a) and (g)). We also note that, for real values of $a$, there is no stability for $b\ge g_2(\alpha)$ (see Figures \ref{fig2}(a), (b), (g), (h), and (n)).\\
2. In the region $g_1(\alpha)>b>g_2(\alpha)$ we obtain two disconnected stable regions (for the complex values of $a$) and one unstable region enclosing an interval on the real axis (see Figure \ref{fig2}(b)).\\
3. If $g_{2}(\alpha)>b>g_{3}(\alpha)$, we obtain a simple closed curve $\gamma(t)$ enclosing the stable region (see Figure \ref{fig2}(c)). This is similar to the case of no feedback ($b=0$).\\
4. In a range between $g_{3}(\alpha)>b>g_{4}(\alpha),$ ($0<\alpha<\alpha_{*}$), $\gamma(t)$ generates one stable region and two (complex valued) unstable regions (see Figure \ref{fig2}(d)). \\
5. In the region $g_{4}(\alpha)>b>g_{3}(\alpha)$, ($\alpha_{*}<\alpha<1$), there is one stable and one unstable region formed by $\gamma(t)$. Both of these regions enclose the real axis (see Figure \ref{fig2}(e)). \\
6. In the region $E$, there are three unstable regions and one stable region (see Figure \ref{fig2}(f)).\\
7. In the region $F$, there are three unstable regions and one stable region (see Figure \ref{fig2}(g)).

\section{Stability result for $\tau=2$}
\begin{figure*}[h]
	\centering
	\includegraphics[scale=0.45]{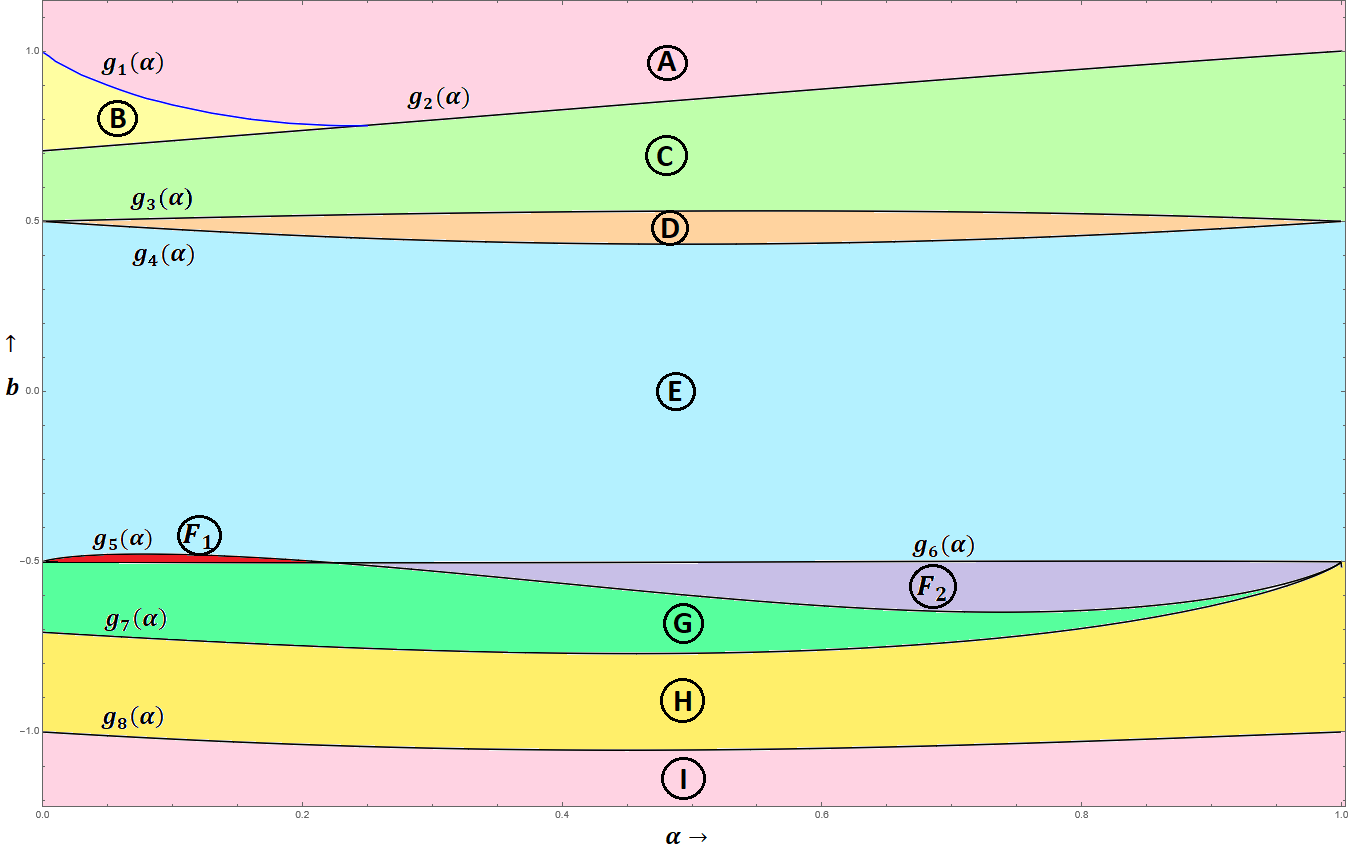}
	\caption{bifurcation regions in $\alpha b$-plane for
		$\tau=2$.}
	\label{fig3}
\end{figure*}

\begin{figure*}[h]
	\centering\subfloat[Region $A$]{%
		\includegraphics[scale=0.22]{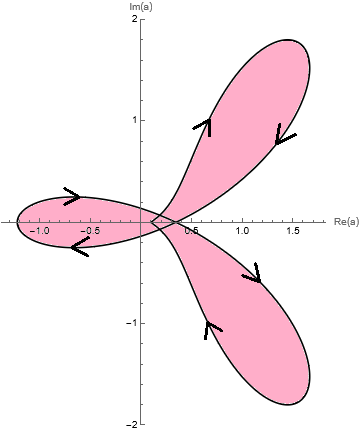}
	}
	\hspace{0.5cm}
	\centering\subfloat[Region $B$]{%
		\includegraphics[scale=0.32]{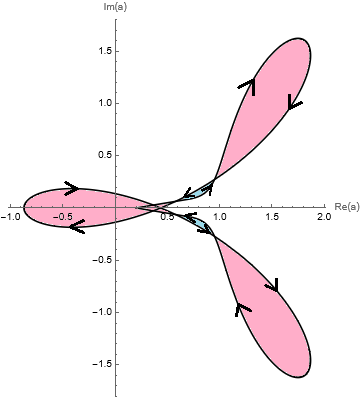}
	}
	\hspace{0.5cm}
	\centering\subfloat[Region $C$]{%
		\includegraphics[scale=0.22]{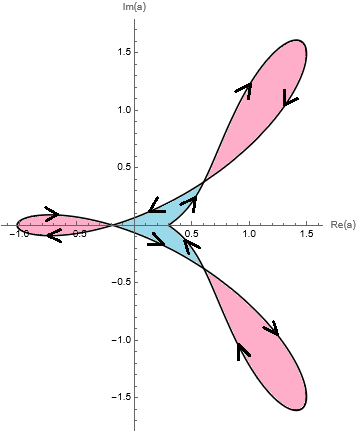}
	}
	\hspace{0.5cm}
	\centering\subfloat[Region $D$]{%
		\includegraphics[scale=0.22]{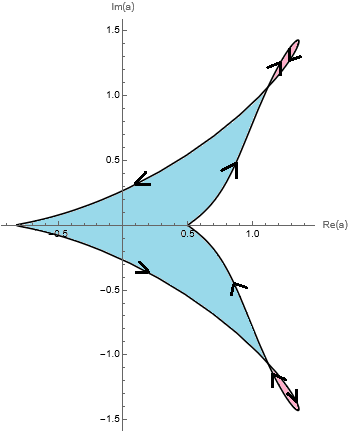}
	}\\
	\centering\subfloat[Region $E$]{%
		\includegraphics[scale=0.22]{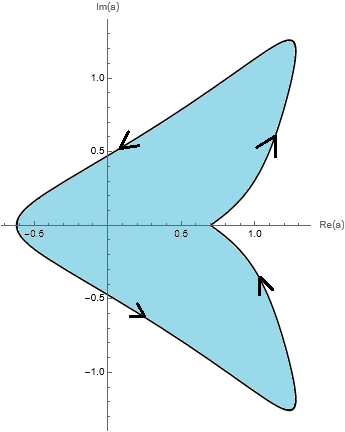}
	}
	\hspace{0.5cm}
	\centering\subfloat[Region $F_{1}$]{%
		\includegraphics[scale=0.22]{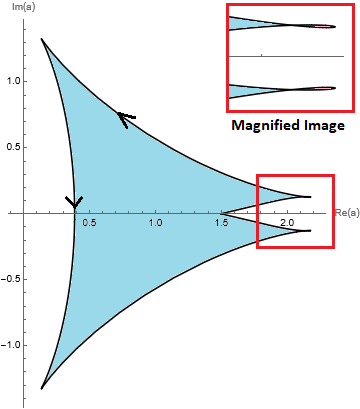}
	}
	\hspace{0.5cm}
	\centering\subfloat[Region $F_{2}$]{%
		\includegraphics[scale=0.22]{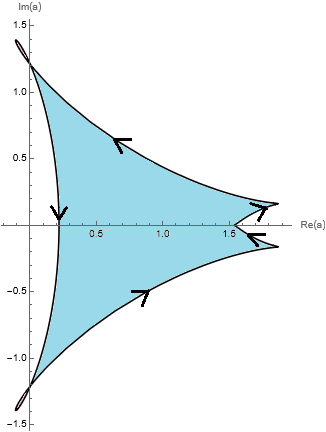}
	}\hspace{0.5cm}
	\centering\subfloat[Region $G$]{%
		\includegraphics[scale=0.22]{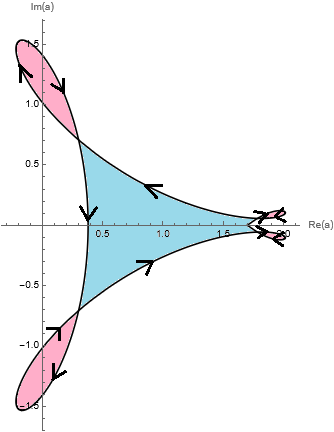}
	}\\
	\centering\subfloat[Region $H$]{%
		\includegraphics[scale=0.22]{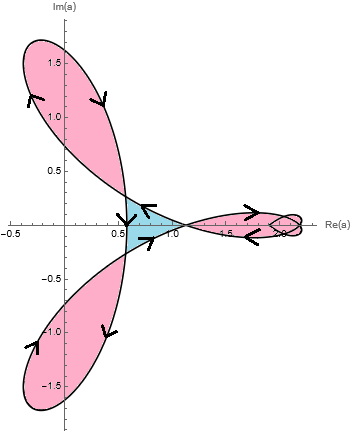}
	}\hspace{0.5cm}
	\centering\subfloat[Region $I$]{%
		\includegraphics[scale=0.22]{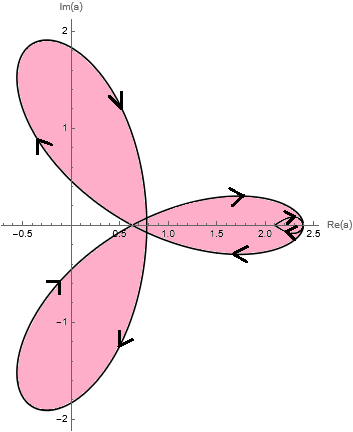}
	}\hspace{0.5cm}
	\centering\subfloat[Curve $b=g_{1}(\alpha)$]{%
		\includegraphics[scale=0.23]{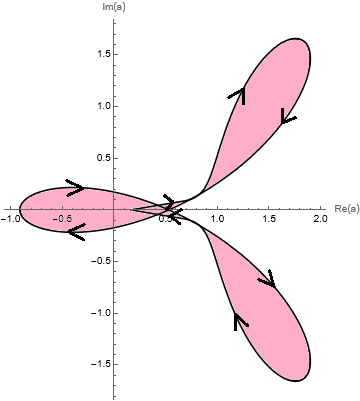}
	}\hspace{0.5cm}
	\centering\subfloat[Curve $b=g_{2}(\alpha)$]{%
		\includegraphics[scale=0.3]{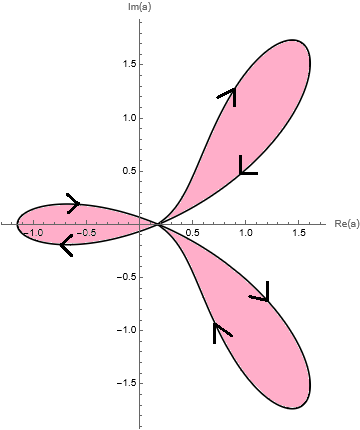}
	}
	\\
	\centering\subfloat[Curve $b=g_{3}(\alpha)$]{%
		\includegraphics[scale=0.3]{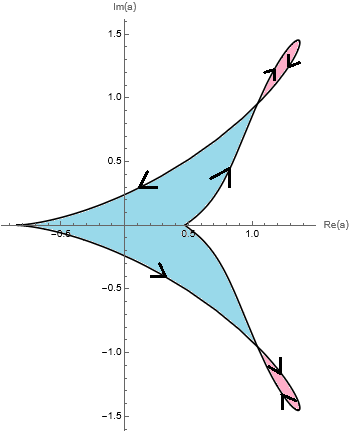}
	}\hspace{0.5cm}
	\centering\subfloat[Curve $b=g_{4}(\alpha)$]{%
		\includegraphics[scale=0.3]{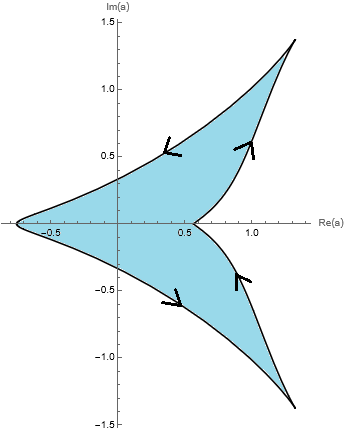}
	}\hspace{0.5cm}
	\centering\subfloat[Curve $b=g_{5}(\alpha)$ ($\alpha<\alpha_{*}$)]{%
		\includegraphics[scale=0.33]{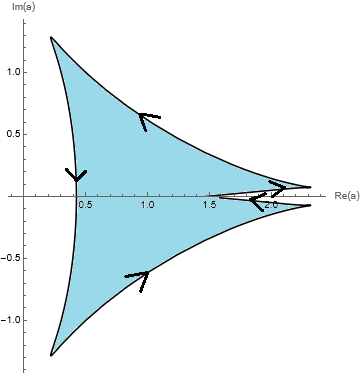}
	}\hspace{0.5cm}
	\centering\subfloat[Curve $b=g_{5}(\alpha)$ ($\alpha>\alpha_{*}$)]{%
		\includegraphics[scale=0.3]{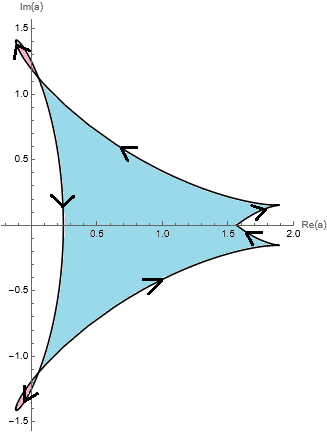}
	}\\
	\centering\subfloat[Curve $b=g_{6}(\alpha)$ ($\alpha<\alpha_{*}$)]{%
		\includegraphics[scale=0.32]{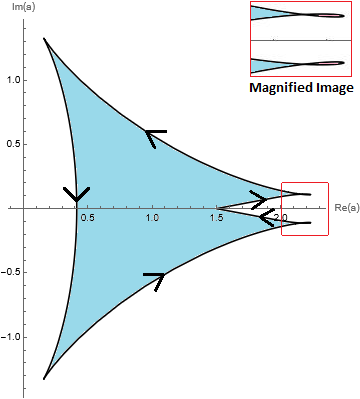}
	}\hspace{0.5cm}
	\centering\subfloat[Curve $b=g_{6}(\alpha)$ ($\alpha>\alpha_{*}$)]{%
		\includegraphics[scale=0.22]{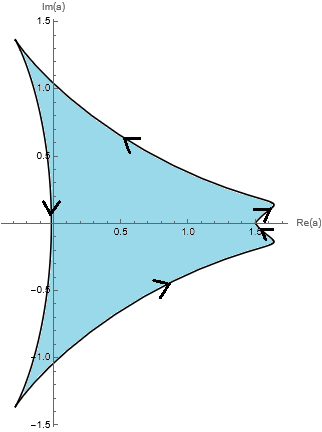}
	}\hspace{0.5cm}
	\centering\subfloat[Curve $b=g_{7}(\alpha)$]{%
		\includegraphics[scale=0.22]{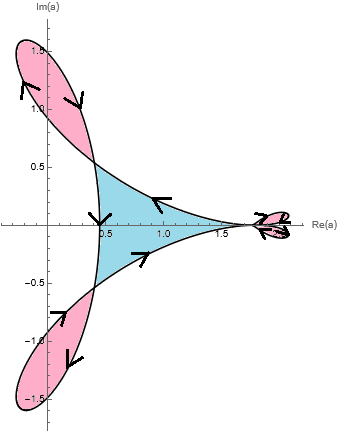}
	}\hspace{0.5cm}
	\centering\subfloat[Curve $b=g_{8}(\alpha)$]{%
		\includegraphics[scale=0.3]{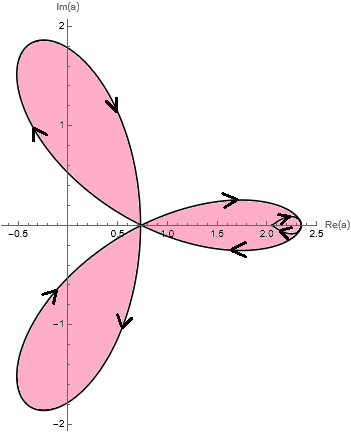}
	}
	\caption{The figures 4(a), 4(b), 4(c), 4(d), 4(e), 4(f), 4(g), 4(h), 4(i), and 4(j), are the representative stability diagrams having $b$ and $\alpha$ values lying in the regions $A, B, C, D, E, F_{1}, F_{2}, G, H,$ and $I$ shown in the figure (\ref{fig3}) respectively. The figures 4(k), 4(l), 4(m), 4(n), 4(o), 4(p), 4(q), 4(r), 4(s) and 4(t) are the representative stability diagrams having $b$ and $\alpha$ values lying on the branches $b=g_{1}(\alpha)$,  $b=g_{2}(\alpha)$, $b=g_{3}(\alpha)$, $b=g_{4}(\alpha)$, $b=g_{5}(\alpha)$, $b=g_{6}(\alpha)$, $b=g_{7}(\alpha)$, and $b=g_{8}(\alpha)$ respectively. Here the blue color indicates stable regions and the pink color indicates unstable regions.}
	\label{fig4}
\end{figure*}

In this section, we provide the bifurcation analysis of the equation (\ref{eqn2}) for $\tau=2$ on the lines of analysis provided for $\tau=1$.
As in the previous case, the curve $\gamma(t)$ has self-intersections for $\tau=2$ also. In this case, there exist distinct points $t_{1},t_{2}\in[0,2\pi]$ such that $\gamma(t_{1})=\gamma(t_{2}).$
We observe that $\gamma(0)=\gamma(t_{1})$, for some $t_{1}$ and $\gamma_{2}(t_{1})=0$. Solving $\gamma_{2}(t_{1})=0$, we get
\begin{equation}
	b=-2^{\alpha}\csc\left(2t_{1}\right)\left(\sin\left(\frac{t_{1}}{2}\right)\right)^{\alpha}\sin\left(\frac{1}{2}\left(t_{1}-t_{1}\alpha+\pi\alpha\right)\right) \label{eqn b3}.
\end{equation}

Substituting this $b$ in $\gamma_{1}(0)=\gamma_{1}(t_{1})$, we get,
\begin{equation*}
	2^{\alpha}\cos\left(\frac{1}{2}\left(-t_{1}\left(-4+\alpha\right)+\pi\alpha\right)\right)\sec\left(t_{1}\right)\left(\sin\left(\frac{t_{1}}{2}\right)\right)^{\alpha}=0.
\end{equation*}
It gives two values of $t_{1}$,
\begin{equation}
	t_{1}=\frac{\pi(-3+\alpha)}{(-4+\alpha)}, \label{g2}
\end{equation}
and
\begin{equation}
	t_{1}=\frac{\pi(-1+\alpha)}{(-4+\alpha)}, \label{g7}
\end{equation}
Using these values of $t_{1}$ in (\ref{eqn b3}) we get the curves
\begin{equation*}
	b=g_{2}(\alpha)=\frac{-2^{\alpha}\left(\cos\left(\frac{\pi}{8-2\alpha}\right)\right)^{\alpha}\cos\left(\frac{\pi}{-4+\alpha}\right)}{\sin\left(\frac{2\pi\left(-3+\alpha\right)}{-4+\alpha}\right)}.
\end{equation*}
and
\begin{equation*}
	b=g_{7}(\alpha)=\frac{2^{\alpha}\left(\cos\left(\frac{3\pi}{8-2\alpha}\right)\right)^{\alpha}\cos\left(\frac{3\pi}{-4+\alpha}\right)}{\cos\left(\frac{3\pi\alpha}{8-2\alpha}\right)}.
\end{equation*}

We observe the cusps in $\gamma(t)$ for four different values of $t$. To find the bifurcation curves for cusps, we solve $\gamma_{1}'(t)=0$ and $\gamma_{2}'(t)=0$. The equation $\gamma_2'(t)=0$ leads to
\begin{equation}
	b=-2^{-2+\alpha}\sec\left(2t\right)\left(\sin\left(\frac{t}{2}\right)\right)^{-1+\alpha}\left(\left(-1+\alpha\right)\sin\left(\frac{1}{2}\left(t+\pi\alpha-t\alpha\right)\right)+\sin\left(\frac{1}{2}\left(3t+\pi\alpha-t\alpha\right)\right)\right), \label{eqn b2}
\end{equation}
Substituting this $b$  in the equation $\gamma_{1}'(t)=0$, we get
\begin{equation*}
	2^{\alpha}\left(\cos\left(\frac{1}{2}\left(-t\left(-7+\alpha\right)+\pi\alpha\right)\right)+\left(-1+\alpha\right)\cos\left(\frac{1}{2}\left(-t\left(-5+\alpha\right)+\pi\alpha\right)\right)\right)\sec\left(2t\right)\left(\sin\left(\frac{t}{2}\right)\right)^{-1+\alpha}=0.
\end{equation*}

Solving this equation for $t$ and substituting back into the equation (\ref{eqn b2}), we get three branches viz. $b=g_{3}(\alpha)$, $b=g_{4}(\alpha)$, $b=g_{5}(\alpha)$ and $b=g_{6}(\alpha)$ as shown in Figure \ref{fig3}. The curves $b=g_{5}(\alpha)$ and $b=g_{6}(\alpha)$ intersect at $\alpha=\alpha_{*}=0.23$.

Similarly, we observed $\gamma(\pi)=\gamma(t_{2})$, for some $t_{2}$ and $\gamma_{2}(t_{2})=0$. Solving $\gamma_{2}(t_{2})=0$, we get
\begin{equation}
	b=-2^{\alpha}\csc\left(2t_{2}\right)\left(\sin\left(\frac{t_{2}}{2}\right)\right)^{\alpha}\sin\left(\frac{1}{2}\left(t_{2}-t_{2}\alpha+\pi\alpha\right)\right). \label{b3}
\end{equation}
Since, $\gamma(\pi)=\gamma(t_{2})$, we solve $\gamma_{1}(\pi)=\gamma_{1}(t_{2})$ and get
\begin{equation}
	2^{\alpha}\left(1+\cos\left(\frac{1}{2}\left(-t_{2}\left(-4+\alpha\right)+\pi\alpha\right)\right)\sec\left(t_{2}\right)\left(\sin\left(\frac{t_{2}}{2}\right)\right)^{\alpha}\right)=0. \label{b3.1}
\end{equation}
The bifurcation curve $b=g_{8}(\alpha)$ is obtained by solving (\ref{b3.1}) and (\ref{b3}) (cf. Fig. \ref{fig3}).

For $\alpha>\alpha_{*}$, if we take $b>g_2(\alpha)$ (region A) then the boundary curve $\gamma(t)$ has negative orientation and the system (\ref{eqn2}) becomes unstable (Fig. \ref{fig4}(a)). However, for smaller values $0<\alpha<\alpha_{*}$ there are two small subregions bounded by  $\gamma(t)$ that are stable (Fig. \ref{fig4}(b)), even after a self-intersection in $\gamma(t)$. These small regions disappear as we increase the value of the parameter $b$ up to the curve $b=g_1(\alpha)$.  
As shown in the Figure \ref{fig3}, the bifurcation curves $b=g_{1}(\alpha)$ and $b=g_{2}(\alpha)$ coalesce at $\alpha=\alpha_{*}$ and  the region bounded by these curves is region  $B$.  

Our observations are summarized below:\\
1. The boundary curve $\gamma(t)$ shows four regions-- none of which are stable-- in the region  $A$ ({\it{i.e.,}} $b > g_{1}(\alpha)$ (Fig. \ref{fig4}(a)).\\
2. In the region $B$ enclosed by $g_1(\alpha)>b>g_2(\alpha)$ which exists only for $\alpha<\alpha_{*}$, we observe the boundary curve $\gamma(t)$ with two disconnected stable regions in the complex plane and five unstable regions (see Figure \ref{fig4}(b)).\\
3. In the region $C$, {\it{i.e.}} $g_{2}(\alpha)>b>g_{3}(\alpha)$, there exist three disconnected unstable regions and a single stable region which encloses the real axis (see Figure \ref{fig4}(c)).\\
4. In the region $D$, {\it{i.e.}} $g_{3}(\alpha)>b>g_{4}(\alpha)$, the boundary curve $\gamma$ has two unstable regions and one stable region. (see Figure \ref{fig4}(d)).\\
5. In the region $E$, {\it{i.e.}} $g_{5}(\alpha)<b<g_{4}(\alpha)$ with $0<\alpha<\alpha_{*}$ and $g_{6}(\alpha)<b<g_{4}(\alpha)$ with $\alpha_{*}<\alpha<1$, the boundary curve $\gamma$ obtained is a simple closed curve with no cusps. The region enclosed by the boundary curve $\gamma$ is stable and the region outside is found to be unstable which is similar to the $b=0$ case (see Figure \ref{fig4}(e)).\\
6. In the region $F_{1}$, {\it{i.e.}} $g_{6}(\alpha)<b<g_{5}(\alpha)$ with $\alpha<\alpha_{*}$, the boundary curve $\gamma$ encloses one stable region and two unstable regions on the right side (see Figure \ref{fig4}(f)).\\
7. In the region $F_{2}$, {\it{i.e.}} $g_{5}(\alpha)<b<g_{6}(\alpha)$ with $\alpha>\alpha_{*}$, the boundary curve $\gamma$ encloses one stable region and two unstable regions on the left side (see Figure \ref{fig4}(g)).\\
8. In the region $G$, {\it{i.e.}} $g_{7}(\alpha)\leq b<g_{6}(\alpha)$  with $0<\alpha<\alpha_{*}$ and $g_{7}(\alpha)\leq b<g_{5}(\alpha)$  with $\alpha_{*}<\alpha<1$, the boundary curve $\gamma$ has four unstable regions and one stable region in the center (see Figure \ref{fig4}(h)).\\
9. In the region $H$, {\it{i.e.}} $g_{8}(\alpha)<b<g_{7}(\alpha)$, the boundary curve $\gamma$ has one stable region and two unstable regions on the left side and the two unstable regions on the right side starts overlapping.\\
10. For the region $I$, {\it{i.e.}} $b<g_8(\alpha)$, there are no stable regions enclosed by $\gamma$.

Now we discuss what happens on the bifurcation curves and a specific point $P$ where $g_5(\alpha)=g_6(\alpha)$.\\
1. The point $P$ where the curves $g_{5}(\alpha)$ and $g_{6}(\alpha)$ intersect {\it{i.e.,}} $(\alpha,b)=(0.23,-0.5)$, the boundary curve $\gamma$ obtained is simple curve with four cusps.\\
2. On the bifurcation curve $b=g_1(\alpha)$, there are six different regions none of which are
stable (Figure \ref{fig4}(k)).\\
3.  On the bifurcation curve $g_{2}(\alpha)$, the $\gamma$ has different number of stable and unstable regions depending on  $\alpha$. If $\alpha<\alpha_{*}$, there are two stable regions and three unstable regions. On the other hand, if $\alpha>\alpha_{*}$, the boundary curve $\gamma$ has three unstable regions(see Figure \ref{fig4}(l)).\\
4. On the bifurcation curve $b=g_{3}(\alpha)$, we observe a cusp on the left side, two unstable and one stable region (see Figure \ref{fig4}(m)).\\
5. On the bifurcation curve $b=g_{4}(\alpha)$ (see Figure \ref{fig4}(n)), and $b=g_5(\alpha)$ (see Figure \ref{fig4}(o) and (p)), there are two cusps on right hand side for the boundary curve $\gamma$.  (see Figure \ref{fig4}(q) and (r)).
On the curve, $g_{5}(\alpha)$, two cusps appear on the right side of the boundary curve $\gamma$. If $\alpha<\alpha_{*}$, the boundary curve $\gamma$ has a simple closed region that is stable (see Figure \ref{fig4}(o)). If $\alpha>\alpha_{*}$, two unstable regions are observed on the left side (see Figure \ref{fig4}(p)).\\
6. We observe two cusps on the left hand side of the boundary for $b=g_6(\alpha)$ (see Figure \ref{fig4}(r)). If $\alpha<\alpha_{*}$, the boundary curve $\gamma$ has two unstable regions on the right side (see Figure \ref{fig4}(q)). But if $\alpha>\alpha_{*}$, the boundary curve $\gamma$ has a simple closed region which is stable (see Figure \ref{fig4}(r)).\\
7. For $b=g_{7}(\alpha)$, the two unstable regions on the right side come closer to forming a cusp (see Figure \ref{fig4}(s)).\\
8. On the bifurcation curve $b= g_{8}(\alpha)$, the stable region enclosed by the boundary curve $\gamma$ disappears and the system becomes completely unstable (see Figure \ref{fig4}(t)).

\section{Nonlinear Maps}
Consider a fractional order map with a delay term 
\begin{equation}
	x(t)=x_0 +\frac{1}{\Gamma(\alpha)}\sum_{j=1}^{t}\frac{\Gamma(t-j+\alpha)}{\Gamma(t-j+1)}[f(x(j-1))+ bx(j-\tau-1)-x(j-1)]. \label{eqn 8}
\end{equation}

If $f(x^*)=(1-b)x^*$, and initial conditions are $x_0=x^*$, $x(t)=x_0=x^*$ for all $t=-\tau\ldots 0$, we obtain $x(t)=x^*$. Thus, $x^*$ is
an equilibrium point. We observed that the stable regions of the linear system (\ref{eqn 7}) are also the stable regions equilibrium of $x^*$ with $a=f'(x^*)$. The numerics show that this is valid for both $x^*=0$ as well as $x^*\ne 0$. We consider only real values of $a$ and $b$.

The system (\ref{eqn2}) having the boundary curve $\gamma(t)$ defined by (\ref{eqn 7}) has stable region given by the curves 
\begin{eqnarray}
	a&=&1-b, \label{a1}\\
	a&=&1-2^{\alpha}-(-1)^{\tau}b \label{a2}
\end{eqnarray}
and the parametric curve 
\begin{equation}
	b(t)=\frac{-2^{\alpha}\left(\sin\left(\frac{t}{2}\right)\right)^{\alpha}\sin\left(\frac{\pi\alpha}{2}+\left(1-\frac{1}{2}\right)t\right)}{\sin\left(\tau t\right)},
	a(t)=1+2^{\alpha}\csc\left(\tau t\right)\left(\sin\left(\frac{t}{2}\right)\right)^{\alpha}\sin\left(\frac{1}{2}\left(-t\left(-2\left(\tau+1\right)+\alpha\right)+\pi\alpha\right)\right) \label{para}
\end{equation} 
where, $t\in[0,2\pi]$ in the $b-a$ plane.

For $\tau=1$, the stable region has 3 corners while for $\tau>1$, the stable region has
four corners. There is a noticeable difference between the stability region of odd and even values of
$\tau$. However, for large values of $\tau$, we observe certain asymptotic stability regions. 

\subsection{Nonlinear maps with fixed point $x^{*}=0$}
\begin{figure*}
	\centering
	\subfloat[Logistic map with $\alpha=0.5$, $\tau=1$.]{%
		\includegraphics[scale=0.3]{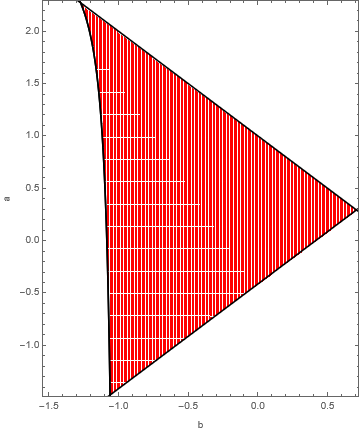}
	}
	\hspace{0.4cm}
	\subfloat[Logistic map with $\alpha=0.75$, $\tau=1$.]{%
		\includegraphics[scale=0.3]{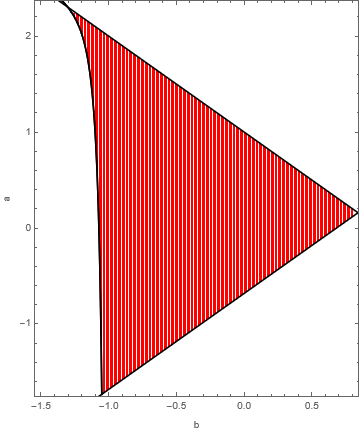}
	}
	\hspace{0.4cm}
	\subfloat[Logistic map with $\alpha=0.5$, $\tau=2$.]{%
		\includegraphics[scale=0.3]{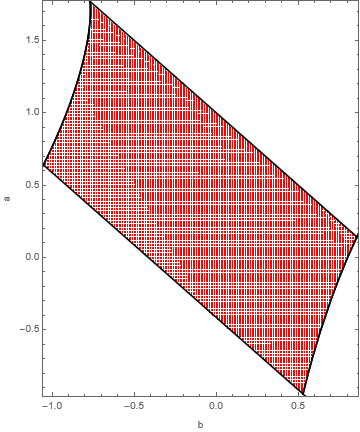}
	}\\
	\subfloat[Cubic map with $\alpha=0.5$, $\tau=1$.]{%
		\includegraphics[scale=0.3]{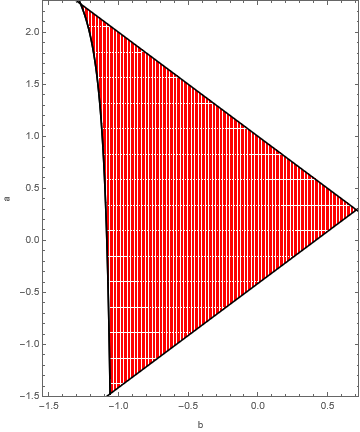}
	}
	\hspace{0.4cm}
	\subfloat[Cubic map with $\alpha=0.75$, $\tau=1$.]{%
		\includegraphics[scale=0.3]{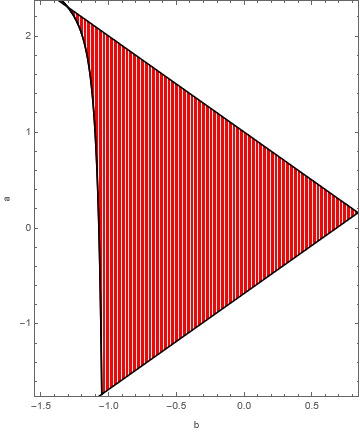}
	}
	\hspace{0.4cm}
	\subfloat[Cubic map with $\alpha=0.5$, $\tau=2$.]{%
		\includegraphics[scale=0.3]{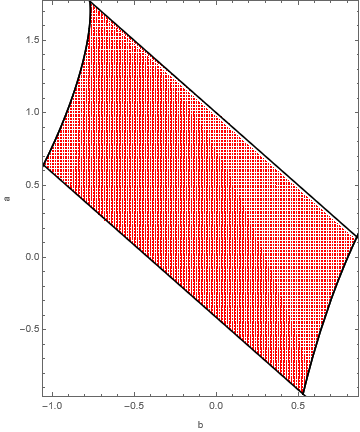}
	}
	\caption{$b$ versus $a$ plot of the systems defined by (\ref{eqn 8}) with logistic and cubic maps. We plot the stable fixed points of the logistic and cubic map }
	\label{fig5}
\end{figure*}

Consider the nonlinear system (\ref{eqn 8}).
Nonlinear maps such as logistic and cubic maps have been studied thoroughly since long and have $0$ as one of the fixed points. For a logistic map, the system (\ref{eqn 8}), can be  explicitly written as:
\begin{equation*}
	x(t)=x_0 +\frac{1}{\Gamma(\alpha)}\sum_{j=1}^{t}\frac{\Gamma(t-j+\alpha)}{\Gamma(t-j+1)}[\lambda x(j-1)(1-x(j-1))+ bx(j-\tau-1)-x(j-1)],
\end{equation*}
where $\lambda$ is the parameter. For a cubic map, the system can be given as:
\begin{equation*}
	x(t)=x_0 +\frac{1}{\Gamma(\alpha)}\sum_{j=1}^{t}\frac{\Gamma(t-j+\alpha)}{\Gamma(t-j+1)}[ \beta x^{3}(j-1)+(1-\beta)x(j-1)+ bx(j-\tau-1)-x(j-1)],
\end{equation*} 
where $\beta$ is the parameter. 

The parameter values  $a=f'(0)$ for the logistic and cubic maps are $\lambda$ and $(1-\beta)$, respectively. We plot the stability region of the system (\ref{eqn 8}) in the $b-a$ plane using the curves (\ref{a1}), (\ref{a2}) and the parametric curve (\ref{para}) for different values of $\alpha$ and $\tau$ for logistic and cubic maps.
Figures \ref{fig5}(a), (b), (d), and (e) show the stability region for  $\tau=1$ and $\alpha=0.5,0.75$ for logistic and cubic maps. The systems are iterated for $10^4$ time-steps where the stable fixed points converge with the precision $\delta=10^{-4}$. From figures \ref{fig5}(a)-(f), it can be observed that the stable region lies within the region enclosed by the $b-a$ curves. This indicates that the range of the stable fixed point can be known for any nonlinear system with $x^*=0$ for given order $\alpha$, time delay $\tau$, and delay coefficient $b$ from the region enclosed by the above boundary curves.
Let us denote the maximum and the minimum values of the parameter $a$ that can be stabilized by strength $b$ as $a_{max}(b)$ and $a_{min}(b)$ respectively. Let $a_{max}=\sup{a_{max}(b)}$ and $a_{min}=\inf{a_{min}(b)}$. For $\tau=1$, the stable range of $a$ is larger for negative values of $b$. We obtain $a_{max}$, when the curve (\ref{a1}) intersects the parametric curve (\ref{para}) and $a_{min}$, when the curve (\ref{a2}) cuts the parametric curve (\ref{para}). Thus we can say that for $\tau=1$, we can stabilize the system (\ref{eqn 8}) for $a\in (a_{min},a_{max}) $ using appropriate value of $b$.

Now, we consider $\tau=2$. In the figures \ref{fig5}(c) and (f), we show the stability regions of logistic and cubic maps for $\tau=2$ and $\alpha=0.5$. If we run simulations for logistic and cubic maps respectively, we find that the entire stability region is filled. We again observe that the theoretical stability regions described by the curves (\ref{a1}), (\ref{a2}), and (\ref{para}) exactly match with those obtained by simulations, for both these maps. Thus the above stability criterion works for nonlinear maps for a fixed point $x^*=0$.

\subsection{Nonlinear maps with nonzero fixed points}
\begin{figure*}
	\centering
	\subfloat[]{%
		\centering\includegraphics[scale=0.23]{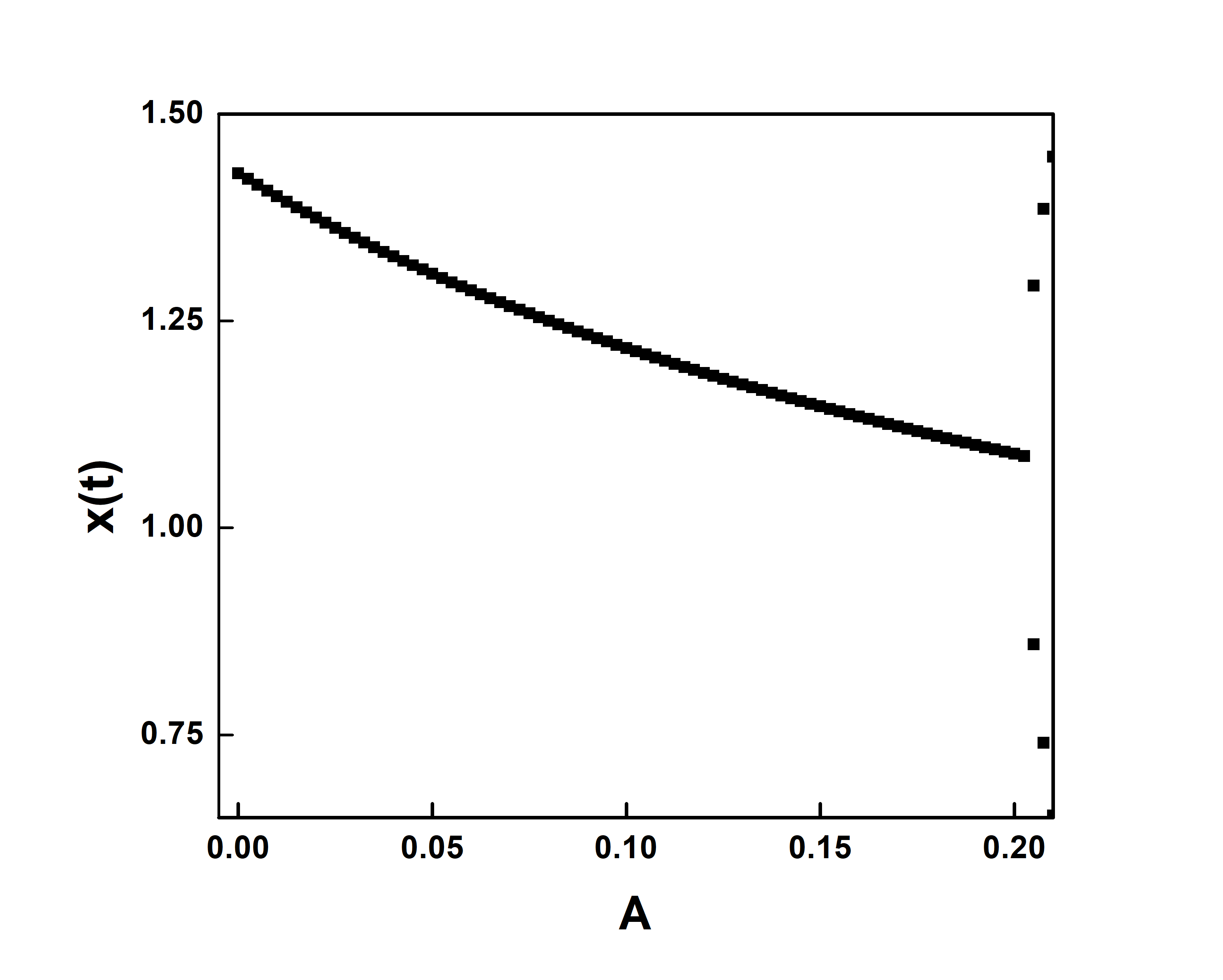}
	}
	\subfloat[]{%
		\centering\includegraphics[scale=0.23]{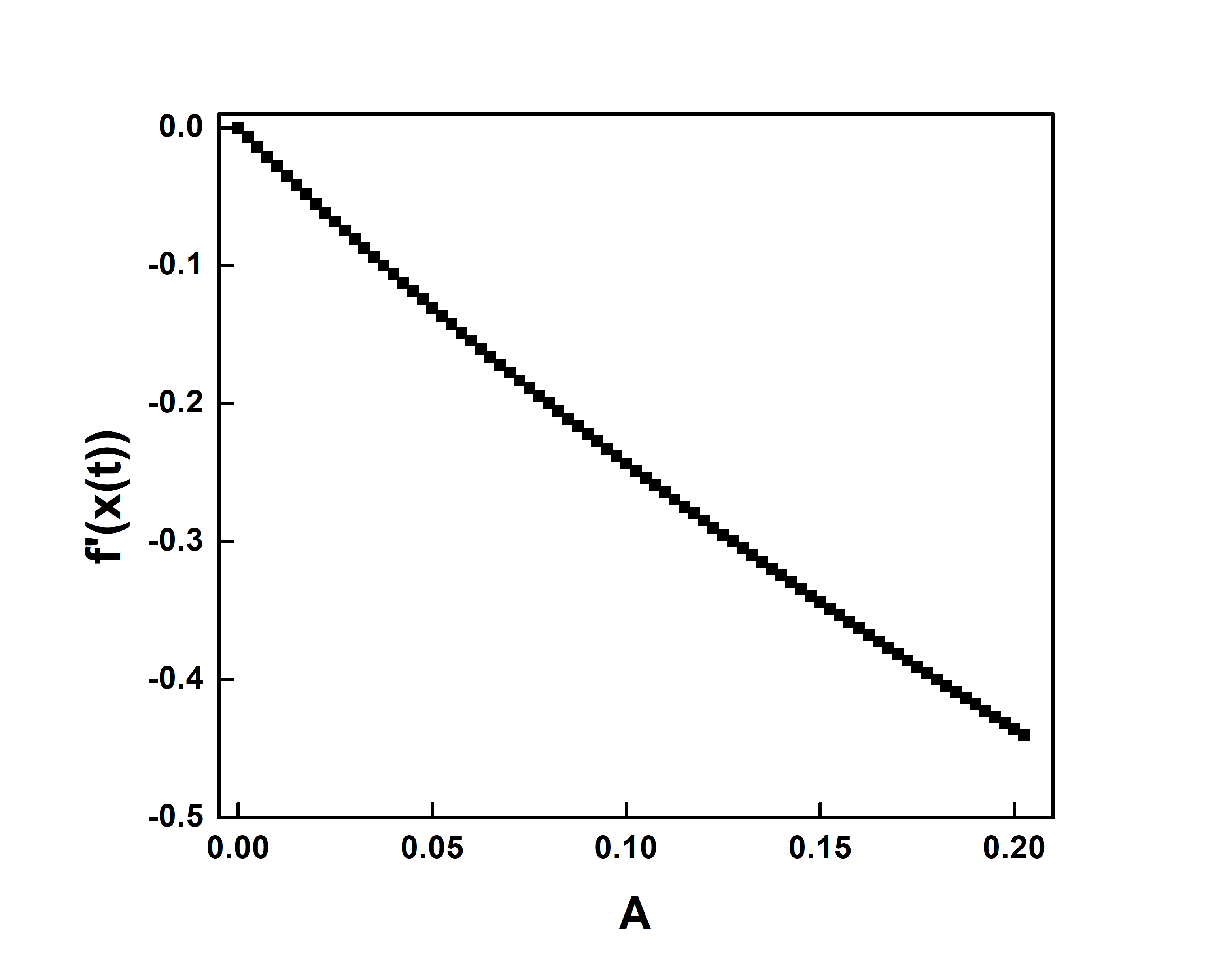}
	}
	\subfloat[]{%
		\centering\includegraphics[scale=0.3]{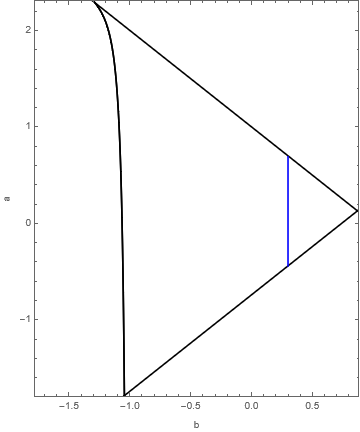}
	}\\
	\subfloat[]{%
		\centering\includegraphics[scale=0.23]{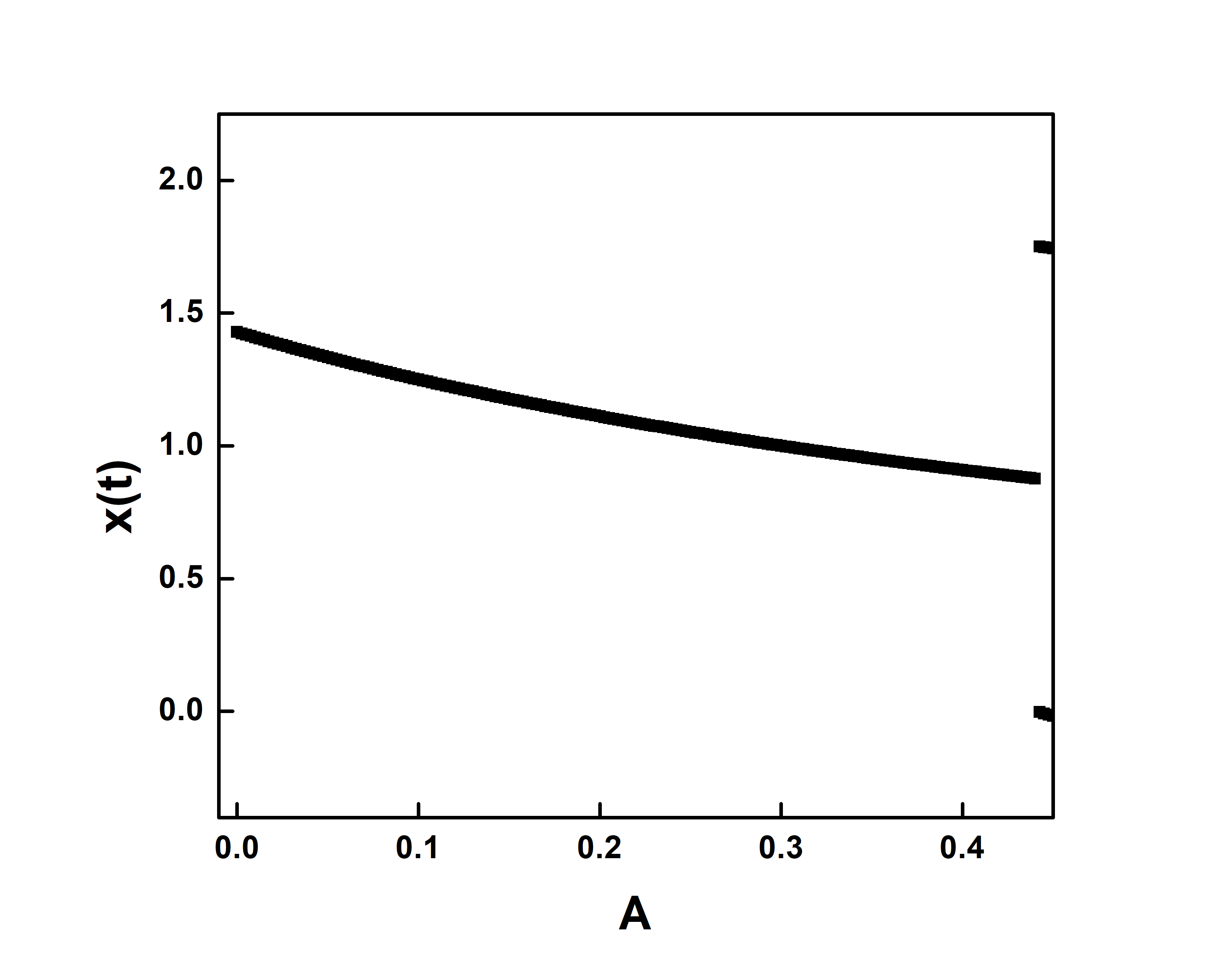}
	}
	\subfloat[]{%
		\centering\includegraphics[scale=0.23]{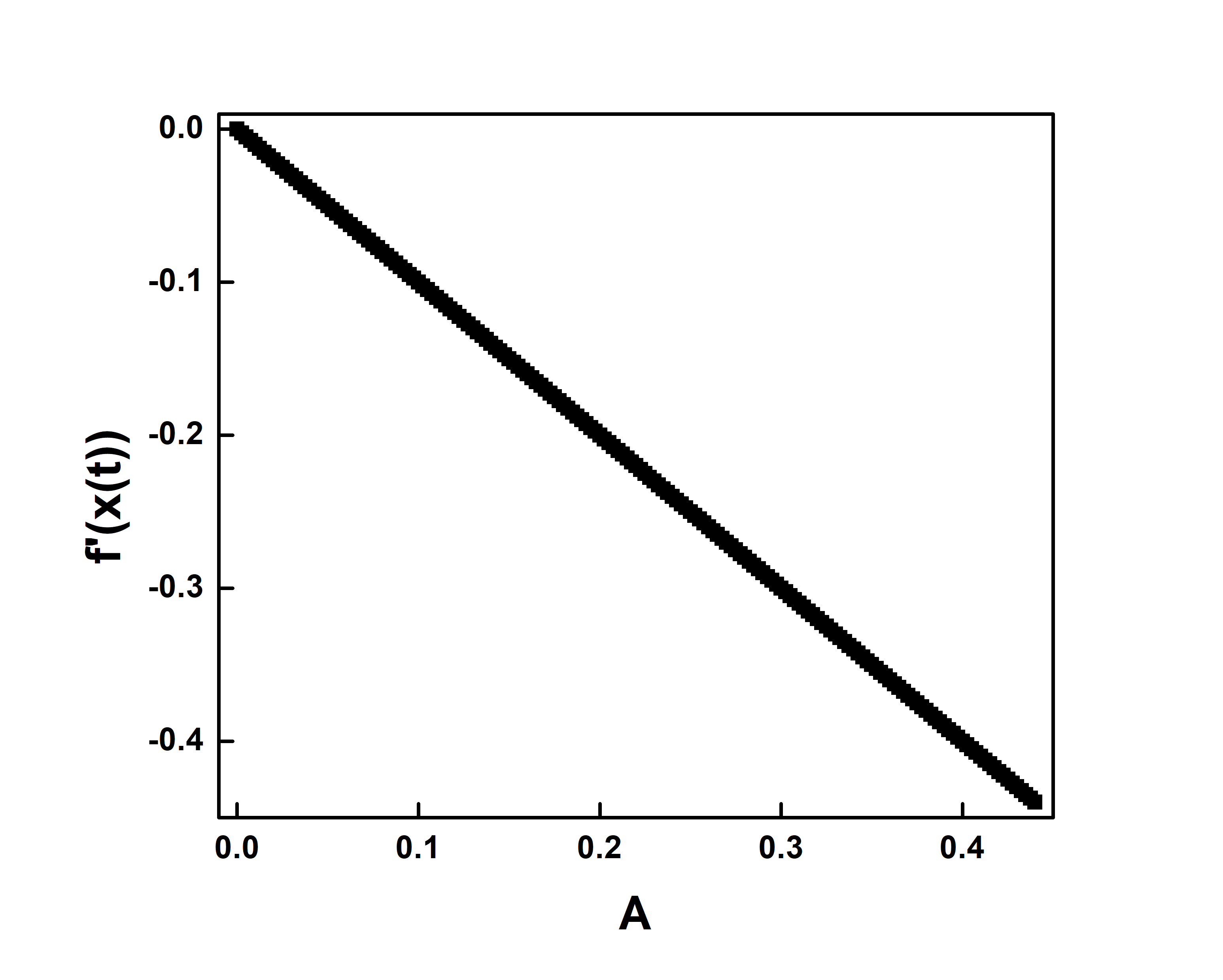}
	}
	\subfloat[]{%
		\centering\includegraphics[scale=0.3]{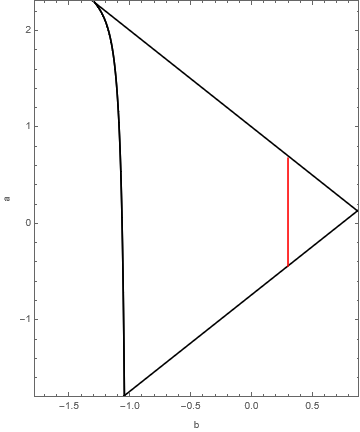}
	}
	\caption{Nonlinear systems defined by (\ref{henon1d}) {\it{i.e.,}} H{\'e}non map and (\ref{lozi1d}) {\it{i.e.,}} Lozi map are iterated for time steps $T=10^{4}$ and the last 120 points are plotted. Here, $\alpha=0.8$, $\tau=1$ and $b=0.3$. Figures (a) and (d) show the range of the stable fixed points of the nonlinear systems (\ref{henon1d}) and (\ref{lozi1d}) respectively. Figures (b) and (e) show the plot for $f'(x(t))$ versus parameter $A$ for the systems (\ref{henon1d}) and (\ref{lozi1d}) respectively. Figures (c) and (f) show the $b$ versus $a$ boundary curve for the systems (\ref{henon1d}) and (\ref{lozi1d}) with $\alpha=0.8$ and $\tau=1$. In both the systems, the line at $b=0.3$ denotes that the stable fixed points of these systems lie within the stable boundary curves defined by (\ref{a1}), (\ref{a2}) and (\ref{para}).}
	\label{figd}
\end{figure*}

Now, we analyze the stability of the map (\ref{eqn 8}) with nonzero fixed points. H{\'e}non and Lozi maps are two such maps with nonzero fixed points. Both of these maps are two-dimensional nonlinear maps. H{\'e}non map is a two-dimensional dissipative quadratic map with coupled equation \cite{henon}, which is defined as 
\begin{eqnarray}
	x_{n+1}=1+y_{n}-Ax^{2}_{n}; \;
	y_{n+1}=Bx_{n}. \label{henon-1}
\end{eqnarray}
Lozi map is a piecewise linear map \cite{lozi}, which is similar to H{\'e}non map where the quadratic term $x^{2}_{n}$ is replaced with $|x_{n}|$. It is defined as,
\begin{eqnarray}
	x_{n+1}=1+y_{n}-A|x_{n}|; \;
	y_{n+1}=Bx_{n}. \label{lozi-1}
\end{eqnarray} 
The equations (\ref{henon-1}) and (\ref{lozi-1}) can also be equivalently written as
\begin{equation}
	x_{n+1}=1-Ax_{n}^2+Bx_{n-1}, \label{h1}
\end{equation}
and 
\begin{equation}
	x_{n+1}=1-A|x_{n}|+Bx_{n-1}, \label{l1}
\end{equation}
respectively.
From the equations (\ref{h1}) and (\ref{l1}), it is observed that both equations can be identified as a delay difference equation with delay one. The fractional order difference equation for (\ref{h1}) and (\ref{l1}) is given below:
\begin{eqnarray}
	x(t)=x(0)+\frac{1}{\Gamma(\alpha)}\sum_{j=1}^{t}\frac{\Gamma(t-j+\alpha)}{\Gamma(t-j+1)}  \left[1-Ax^{2}(j-1)+Bx(j-2)-x(j-1)\right], \label{henon1d}
\end{eqnarray}
and 
\begin{eqnarray}
	x(t)=x(0)+\frac{1}{\Gamma(\alpha)}\sum_{j=1}^{t}\frac{\Gamma(t-j+\alpha)}{\Gamma(t-j+1)} \left[1-A\vert x(j-1)\vert+Bx(j-2)-x(j-1)\right], \label{lozi1d}
\end{eqnarray}
respectively.
As expected, the equations (\ref{henon1d}) and (\ref{lozi1d}) resemble (\ref{eqn 8}) with $\tau=1$. Thus, the stability of the fixed points of these maps can be analyzed similarly. For H{\'e}non map, we identify $f(x)=1-A x^2$ and $b=B$ in equation (\ref{henon1d}). For Lozi map, we identify, $f(x)=1-A\vert x\vert$ and $b=B$ in equation (\ref{lozi1d}). For $x^*$ such that $f(x^*)+bx^*=x^*$, if we take $x(0)=x(-k)=x^*$ for $k=1\ldots \tau$, it is an equilibrium point. The question is if the stability boundary of a fixed point of the above system is given by $a=f'(x^*)$ and $b$ values and if the linearization holds. We plot the bifurcation diagrams for (\ref{henon1d}) and (\ref{lozi1d}) with $\alpha=0.8$, $B=0.3$. The systems are iterated for $10^4$ time-steps where the stable fixed points converge with the precision $\delta=10^{-5}$. For $b=B=0.3$, the stability limits of values of $a_{max}(b)=0.7$ and $a_{min}(b)=-0.441$ for $\tau=1$. It is observed that the range of fixed point of the systems (\ref{henon1d}) and (\ref{lozi1d}) lies within the stable region described for $\tau=1$ and $\alpha=0.8$ (see Figure \ref{figd}) for initial conditions close to the fixed point. Thus we conjecture that having $(b, a)$ values in the stable domain in the $b-a$ plane is a necessary and sufficient condition for the stability of nonzero fixed points of nonlinear maps as well.

\section{Discussion and Conclusion}

Delay difference equations as well as fractional difference equations have memory. Delay difference equations have finite memory. The difference equation with delay $\tau$ is equivalent to the $\tau+1$ dimensional difference equation. The stability conditions for such equations are similar to that for $\tau+1$ dimensional difference equations. These equations are useful in several physical systems. The fractional difference equations, on the other hand, have a long-term memory. We find that the stability properties of fractional order systems with delay are very different from $\tau+1$ dimensional fractional difference equations. For $\tau+1$ dimensional fractional difference equations, the stability is guaranteed if all eigenvalues of corresponding Jacobian lie in the simply connected stable region \cite{bhalekar2022stability}. The system shows a complex set of bifurcations which become more complicated for higher $\tau$. 

In this work, we have studied linear fractional difference equations with a delay term. We derived the boundary curve (\ref{eqn 7}) in a complex plane for the system (\ref{eqn2}). The detailed stability analysis is given for the cases $\tau=1$ and $\tau=2$ and an extremely complex set of bifurcations is obtained.  

We further extended this investigation for the nonlinear map with the fixed point $x^*=0$ as well as $x^*\ne 0$. We found that by identifying $a=f'(x^*)$ same stability 
conditions are useful in analyzing the stability of delayed nonlinear maps.

\section{Acknowledgement}
PMG and DDJ thank DST-SERB for financial assistance (Ref. CRG/2020/003993).

\bibliographystyle{plain}

\bibliography{newref}

\end{document}